\documentclass[reqno]{amsart}

\input xy
\xyoption{all}
\usepackage{accents}
\usepackage{epsfig}
\usepackage{color}
\usepackage{amsthm}
\usepackage{amssymb}
\usepackage{amsmath}
\usepackage{amscd}
\usepackage{amsopn}
\usepackage{graphicx}
\usepackage{tikz-cd} 

\usepackage{xspace}

\usepackage{url}
\usepackage{enumitem, hyperref}\hypersetup{colorlinks}

\usepackage{color} 

\definecolor{darkred}{rgb}{1,0,0} 
\definecolor{darkgreen}{rgb}{0,0.8,0}
\definecolor{darkblue}{rgb}{0,0,1}

\hypersetup{colorlinks,
linkcolor=darkblue,
filecolor=darkgreen,
urlcolor=darkred,
citecolor=darkgreen}

\makeatletter
\def\reflb#1#2{\begingroup
    #2%
    \def\@currentlabel{#2}%
    \phantomsection\label{#1}\endgroup
}
\makeatother

\numberwithin{equation}{section}
\newtheorem {Theorem}{Theorem}
\numberwithin{Theorem}{section}

\newtheorem {Corollary}[Theorem]{Corollary}
\theoremstyle{definition}

\theoremstyle{remark}
\newtheorem{Remark}[Theorem]{Remark}

\newcommand{\CS}{{\mathcal S}}

\newcommand{\Aa}{{\mathcal A}}

\newcommand{\Ss}{{\mathcal S}}

\def    \nat    {{\natural}}

\def    \R      {{\mathbb R}}
\def    \reals      {{\mathbb R}}
\def    \Z      {{\mathbb Z}}
\def    \N      {{\mathbb N}}

\def    \12    {{\frac{1}{2}}}

\def    \HF     {\operatorname{HF}}

\def    \H     {\operatorname{H}}

\def    \CF     {\operatorname{CF}}

\def    \FA     {\operatorname{FA}}

\def    \bx     {\bar{x}}

\def    \MUCZ  {\operatorname{\mu_{\scriptscriptstyle{CZ}}}}

\begin{document}


\setlength{\smallskipamount}{6pt}
\setlength{\medskipamount}{10pt}
\setlength{\bigskipamount}{16pt}





\title[Conley conjecture and local Floer homology]{Conley conjecture and local Floer homology}

\author[Erman \c C\. inel\. i]{Erman \c C\. inel\. i}

\address{Department of Mathematics \\
UC Santa Cruz \\ 
Santa Cruz, CA 95064 \\
USA}
 \email{scineli@ucsc.edu}



\date{\today} 


\begin{abstract} 
In this paper we connect algebraic properties of the pair-of-pants product in local Floer homology and Hamiltonian dynamics. We show that for an isolated periodic orbit the product is non-uniformly nilpotent and use this fact to give a simple proof of the Conley conjecture for closed manifolds with aspherical symplectic form. More precisely, we prove that on a closed symplectic manifold the mean action spectrum of a Hamiltonian diffeomorphism with isolated periodic orbits is infinite. 
\end{abstract}

\maketitle


\tableofcontents

\section{Introduction and main results}
\label{sec:intro+results}

\subsection{Introduction}
\label{sec:intro} 
In this paper we study the pair-of-pants product in local Floer homology and use its properties to give a simple proof of the Conley conjecture for closed manifolds with aspherical symplectic form. We prove that the product in the local Floer homology is non-uniformly nilpotent unless the periodic orbit is a symplectically degenerate maximum. Then we utilize this fact to prove that on a closed symplectic manifold the mean action spectrum of a Hamiltonian diffeomorphism with isolated periodic orbits is infinite. 

To state the results in more detail; recall that for a broad class of symplectic manifolds, every Hamiltonian diffeomorphism has infinitely many simple periodic orbits. Such existence results are usually referred as the Conley conjecture. The example of an irrational rotation of $M=S^2$ shows that the conjecture does not hold unconditionally, or the condition $\omega\vert_{\pi_2(M)}=0$ cannot be dropped completely. The following fact, proven in \cite{GG:conley}, encompasses all known cases of the conjecture: When a closed symplectic manifold $(M,\omega)$ admits a Hamiltonian diffeomorphism with finitely many periodic orbits, there
is a class $A\in\pi_2(M)$ with $\omega(A)>0$ and
$\left<c_1(TM),A\right>>0$. The proof is a formal consequence of previously known cases combined with the aspherical case established in \cite{GG:conley}. Ginzburg and G\" urel prove the aspherical case by constructing a strictly decreasing sequence of mean action values. In this paper, following Remark 4.4 in \cite{GG:conley}, we show that the mean action spectrum is infinite by using a vanishing property of the pair-of-pants product in local Floer homology. 

Here is a brief outline of the argument. Let $\varphi_H$ be a Hamiltonian diffeomorhism of a closed symplectic manifold generated by a one-periodic Hamiltonian $H$. We first show that the local Floer algebra $\FA(H,x) = \bigoplus_{k=1}^{\infty} \HF_*(H^{\nat k}, x^k)$ of a one-periodic orbit $x$ is non-uniformly nilpotent, if $x^k$ is isolated and not a symplectically degenerate maximum for all $k\in \N$ (cf. Prop. 5.3 in \cite{GG:gap}). Then arguing by contradiction, assuming in addition that the mean action spectrum of $\varphi_H$ is finite, we show that the total Floer algebra $\FA(H)=\bigoplus_{k=1}^{\infty} \HF_*(H^{\nat k})$ is nilpotent, which is impossible.

\subsection{Main results}
\label{sec:results} Let us now state the main theorems. The conventions and basic definitions are reviewed in Section \ref{sec:prelim}. In what follows a  ``periodic orbit of a Hamiltonian diffeomorphism" means a  ``contractible periodic orbit of the time-dependent flow generated by a Hamiltonian".  

\begin{Theorem}
\label{thm:1}
Let $\varphi_H$ be a Hamiltonian diffeomorphism of a closed symplectic manifold $(M, \omega)$ generated by a one-periodic Hamiltonian $H$. Assume that periodic orbits of $\varphi_H$ are isolated. Then the mean action spectrum of $H$ is infinite. 
\end{Theorem}

When $\omega$ is not aspherical, a single periodic orbit gives rise to an infinite action spectrum. In that case the assertion is trivial, since every Hamiltonian diffeomorphism of a closed manifold has a periodic orbit. On the other hand, when $\omega$ is aspherical, Theorem \ref{thm:1} implies existence of infinitely many simple periodic orbits. In other words, Theorem \ref{thm:1} implies that the Conley conjecture holds  for closed symplectic manifolds $(M, \omega)$ with aspherical $\omega$.

\begin{Corollary}
\label{cor}
A Hamiltonian diffeomorphism of a closed symplectic manifold $(M, \omega)$ with aspherical $\omega$ has infinitely many simple periodic orbits.
\end{Corollary}

This is a Lusternik-Schnirelmann type result in the sense that the lower bound for critical points is established by bounding critical values. 
Theorem \ref{thm:1} is proved in Section \ref{sec:pf1}. The proof relies on the following vanishing property of the pair-of-pants product.

Let $M$ be a symplectic manifold and $x$ be a one-periodic orbit of a Hamiltonian $H \colon S^1 \times M \rightarrow \reals$. Denote by $x^k$ the $k$th iteration of $x$. The local Floer homology $\HF_*(H^{\nat k}, x^k)$ is defined whenever $x^k$ is isolated. The pair-of-pants product turns the direct sum 
\[
\FA(H,x)=\bigoplus_{k=1}^{\infty} \HF_*(H^{\nat k}, x^k)
\] 
into a graded algebra, which we call the local Floer algebra. We say that a graded algebra $\bigoplus_{k=1}^{\infty} \mathcal{A}_k$ is \emph{non-uniformly nilpotent}, if for all $N \in \N$ there is $m\in\N$ such that all $m$-fold products vanish when restricted to $\bigoplus_{k=1}^N \mathcal{A}_k$. The next result (cf. Prop. 5.3 in \cite{GG:gap}) is proved in Section \ref{sec:pf2}. See Section \ref{sec:floer} for the definition of symplectically degenerate maximum.  

\begin{Theorem}
\label{thm:2}
Let $M$ be a symplectic manifold and let $x$ be a one-periodic orbit of a Hamiltonian $H \colon S^1 \times M \rightarrow \reals$. Assume that $x^k$ is isolated and not a symplectically degenerate maximum for all $k \in  \N$. Then the local Floer algebra $\FA(H,x)$ is non-uniformly nilpotent. 
\end{Theorem}

\medskip

\noindent{\bf Acknowledgements.} The author is deeply grateful to Ba\c sak G\" urel and Viktor Ginzburg for their numerous valuable remarks and suggestions. The author also thank the anonymous referee for critical remarks and helpful comments. A part of this work was carried during \emph{Hamiltonian and Reeb dynamics: New methods and Applications} workshop at Lorentz Center (Leiden, Netherlands); and supported by UCSC. The author would like to thank these institutions for their support. 


\section{Preliminaries}
\label{sec:prelim}

\subsection{Conventions and basic definitions} 
\label{sec:notation}
Let $(M,\omega)$ be a closed symplectic manifold and $H$ be a one-periodic in time Hamiltonian on $(M, \omega)$, i.e., $H \colon S^1 \times M \rightarrow \reals$, where $S^1=\reals/ \Z$. The Hamiltonian vector field $X_H$ of $H$ is defined by $i_{X_H}\omega= -dH$. The time-one map of the time-dependent flow of $X_H$ is denoted by $\varphi_H$. Such time-one maps are referred as \emph{Hamiltonian diffeomorphisms}.

A \emph{capping} of a contractible loop $x \colon S^1 \rightarrow M$ is a map $A \colon D^2 \rightarrow M$ such that $A\vert_{S^1} =x$. The action of a Hamiltonian $H$ on a capped closed curve $\bx=(x,A)$ is
\[
\Aa_H(\bx)=-\int_A \omega + \int_{S^1} H(t,x(t)) \, dt. 
\]
The critical points of $\Aa_H$ on the space of capped closed curves are exactly the capped one-periodic orbits of $X_H$. The set of critical values of $\Aa_H$  is called the \emph{action spectrum} $\Ss (H)$ of $H$ (or of $\varphi_H$). These definitions extend to Hamiltonians of any period in an obvious way.

For $k\in \N$, the  $k$th iteration of $H$, by which we simply mean $H$ treated as $k$-periodic, is denoted by $H^{\nat k}$. With this notation, the \emph{mean action spectrum} $\widehat{\CS}(H)$ of $H$ (or of $\varphi_H$) is defined as the union of the normalized spectra $\CS\big(H^{\nat k}\big)/k$.  Note that the action functional is homogeneous with respect to iteration: 
\[
\Aa_{H^{\nat k}}(\bx^k)=k\Aa_H(\bx)
\]
where $\bx^k$ is the $k$th iteration of the capped orbit $\bx$. Furthermore, attaching a sphere $A \in \pi_2(M)$ to  $\bx^k$ changes $\Aa_{H^{\nat k}}(\bx^k)$ by $\omega(A)$. So $\widehat{\CS}(H)$ is always infinite when $\omega$ does not vanish on $\pi_2(M)$, and when $\omega$ is \emph{aspherical} (i.e. $\omega\vert_{\pi_2(M)}=0$) Theorem \ref{thm:1} implies the existence of infinitely many \emph{simple} (i.e., not iterated) periodic orbits. 

A periodic orbit $x$ of $H$ is called \emph{non-degenerate} if the linearized return map $d\varphi_H \colon T_{x(0)}M \rightarrow  T_{x(0)}M$ has no eigenvalues equal to one. The \emph{Conley-Zehnder index} $\MUCZ (\bx) \in \Z$ of a non-degenerate capped orbit $\bx$ is defined, up to a sign, as in \cite{Sa,SZ}. In this paper, $\MUCZ$ is normalized so that $\MUCZ(\bx)=n$ when $x$ is a maximum, with trivial capping, of an autonomous Hamiltonian with small Hessian. The \emph{mean index} $\Delta (\bx) \in \reals$ is defined even when $\bx$ is degenerate and depends continuously on $H$ and $\bx$, see \cite{Lo,SZ}. Furthermore, it satisfies 
$\big|\Delta(\bx)-\MUCZ(\bx)\big|\leq n$, and it is homogeneous with respect to iteration:
\[
\Delta\big(\bx^k\big)=k\Delta(\bx).
\]

\subsection{Floer homology} In this section, we recall some properties of local Floer homology and pair-of-pants product, then define Floer algebra. See \cite{FO,GG:gaps,HS,MS,Sa,SZ} for a detailed account on (filtered) Floer homology, and \cite{AS,MS,PSS} for pair-of-pants product. 
  
\subsubsection{Local Floer homology} 
\label{sec:floer}  Let $\bx$ be an isolated capped periodic orbit of a Hamiltonian $H$. The local Floer homology $\HF_*(H,\bx)$ of $\bx$ is defined as in \cite{Gi:CC,GG:gaps,GG:gap}.     
The capping is only used to fix a trivialization of $TM\vert_x$ and hence to give a $\Z$-grading to $\HF_*(H,\bx)$. Recapping a capped orbit shifts Conley-Zehnder index by an even integer. So $\Z_2$-graded homology $\HF_*(H,x)$ is defined without fixing a trivialization. 
 
The \emph{support} of $\HF_*(H,\bx)$ is the collection of integers $m$ such that $\HF_m(H,\bx) \neq 0$. A capped periodic orbit $\bx$ is called a \emph{symplectically degenerate maximum (SDM)} if $\Delta(\bx)$ is an integer and $\HF_{\Delta(\bx)+n}(H,\bx) \neq 0$, this property is independent of the capping. The mean action spectrum of a Hamiltonian diffeomorphism $\varphi_H$ with an SDM orbit is infinite, see \cite{Gi:CC} for details. If $\bx$ is not an SDM, then the support of $\HF_*(H,\bx)$  is contained in the half-open interval $[ \Delta(\bx)-n, \Delta(\bx)+n  )$. 

We will use the properties of $\Z$-grading by Conley-Zehnder index in the proof of Theorem \ref{thm:2}. When proving Theorem \ref{thm:1}, for the sake of simplicity, we work with $\Z_2$-graded homology.

\subsubsection{Pair-of-pants product} 
\label{sec:product}
Let $(M,\omega)$ be a closed symplectic manifold with \emph{rational} $\omega$, i.e., $\omega \vert_{\pi_2(M)}$ is discrete. The filtered Floer homology on $(M, \omega)$ carries the so-called \emph{pair-of-pants} product; see, e.g., \cite{AS}.
 On the level of complexes, this product is a map
\[
\CF^{(a,\, b)}_m(H)\otimes \CF^{(c,\, d)}_l(K)
\to \CF^{(\max\{ a+d,\, b+c \},\, b+d)}_{m+l-n}(H\nat K)
\]
giving rise on the level of homology to an associative,
commutative product
\[
\HF^{(a,\, b)}_m(H)\otimes \HF^{(b,\, d)}_l(K)
\to \HF^{(\max\{ a+d,\, b+c \},\, b+d)}_{m+l-n}(H\nat K)
\]
where $H$, $K$ are one-periodic Hamiltonians that satisfy $H(0, \cdot)= K(0, \cdot)$ with all the time derivatives and $H\sharp K$ is the two-periodic Hamiltonian given by $H(t,\cdot)$ for $t\in [0,1]$ and $K(t-1,\cdot)$ for $t\in [1,2]$. Here $H$, $K$ could be Hamiltonians of any period, then  $H\sharp K$  would have the sum of the periods. The product turns the direct sum of the filtered Floer homology
\[
\FA^{(-\infty,\, b)}(H):=\bigoplus_{k =1}^{\infty} \HF^{(-\infty, \, kb)}_*\big(H^{\nat k}\big) 
\]
into an associative and graded-commutative non-unital algebra, which we call the \emph{filtered Floer algebra}. The \emph{local Floer algebra} 
\[
\FA(H,x):=\bigoplus_{k=1}^{\infty} \HF_*(H^{\nat k}, x^k)
\] 
and the \emph{total Floer algebra} $\FA (H)$ are defined in a similar way. In all three, we work with $\Z_2$-graded homology when $c_1(TM)\neq 0$.  In the local case, similar to how local Floer homology is defined, the product is given by counting pair-of-pants curves between the orbits contained (as fixed points) in a small neighborhood of $x$ (as a fixed point). At the end of this section, we justify why this restricted product is a chain map. Note that we only use the rationality assumption for $\omega$ to define filtered Floer homology; in the case of local and total algebra the assumption is not needed.

In Section \ref{sec:pf2} when proving nilpotency of local Floer algebra we use properties of $\Z$-grading by Conley-Zehnder index, in other words we work with the $\Z$-graded local algebra $\FA(H,\bx)$ of a capped orbit $\bx$. Here to cap iterated orbits $x^k$, we iterate the capped orbit $\bx$.  Once we prove the nilpotency result, we will forget $\Z$-grading and go back to $\Z_2$-graded algebra $\FA(H, x)$ which does not depend on the capping. 

The product structure in total Floer algebra, under the isomorphism with the (dual of) quantum cohomology at each summand, agrees with the (dual of) quantum cup-product. Using the unit element of the quantum cup-product at each summand of $\FA (H)$ one can form a non-zero product of any length. We prove Theorem \ref{thm:1} by contradiction: We show that $\FA(H)$ would be nilpotent if the mean action spectrum were finite. 

In Section \ref{sec:pf1}, arguing by contradiction, we will extend the nilpotency result in local Floer algebra to the total algebra under the assumption that the mean action spectrum is finite. Below is the first step towards building total algebra out of local algebras of isolated orbits; one can also use this approach to define pair-of-pants product in the local Floer homology of an isolated orbit under the rationality assumption. Suppose that $\omega$ is aspherical. For an isolated mean action value $c$ one can consider the filtered Floer algebra around $c$:
\[
\FA^{(c-\epsilon, \, c+\epsilon)}(H):=\bigoplus_{k =1}^{\infty} \HF^{(kc-k\epsilon, \, kc+k\epsilon)}_*\big(H^{\nat k}\big).
\]
Here it is essential that the interval $(c-\epsilon, \, c+\epsilon)$ contains a single mean action value. Similar to Floer homology of a single Hamiltonian, as a vector space, the filtered algebra $\FA^{(c-\epsilon, \, c+\epsilon)}(H)$ splits
\[
\FA^{(c-\epsilon, \, c+\epsilon)}(H)= \bigoplus_{\mathcal{A}_{H^{\nat k}}(\bx_i)=kc} \FA(H^{\nat k}, \bx_i)
\]
where $\bx_i$'s are simple (not iterated) capped orbits. This splitting also respects product structure. In particular, one can show that (see \cite{GG:conley}) on the level of complexes there are no pair-of-pants curves between distinct summands of this splitting for sufficiently small non-degenerate perturbations of Hamiltonians. 
Such a pair-of-pants curve would have an a priori energy lowerbound that does not depend on non-degenerate pertubations of the Hamiltonians (see Prop. 2.2 in \cite{GG:conley}). On the other hand, the energy identity for pair-of-pants curves implies that by controlling perturbations, we can keep the energy of all pair-of-pants  curves as small as we want (here we are using the fact that all non-perturbed orbits have the same mean action value). 

Recall that the product in the local algebra $\FA(H,x)$ is given by counting pair-of-pants curves that connects orbits contained in a small neighborhood of $x$. Arguing as in the previous paragraph, by controlling perturbations, one can guarantee that index one pair-of-pants curves cannot brake at an orbit outside of the chosen neighborhood of $x$; which implies that the restricted product is a chain map. Such a curve would have an a priori energy lowerbound that doesn't depend on the perturbations (see Prop. 2.2 in \cite{GG:conley}), and we can keep the energy of all curves as small as we want by controlling perturbations.

\section{Proofs}
\label{sec:pf} 

\subsection{Proof of Theorem \ref{thm:2}}
\label{sec:pf2} 

Let $x$ be a one-periodic orbit of a Hamiltonian $H \colon S^1 \times M \rightarrow \reals$. Assume that $x$ and all of its iterations are isolated, and non-SDM. Choose a capping $\bx =(x, A)$ and iterate it $\bx^k =(x^k,A^k)$. For a fixed  $N \in \mathbb{N}$, consider products of the form $w_1\cdot\dots {\cdot w_r} \in \HF_l(H^{\nat k}, \bx^k)$ with non-zero $w_i \in \HF_{l_i}(H^{\nat k_i}, \bx^{k_i})$; where $k_i \leq N$, $l= \sum l_i -(r-1)  n$ and $k=\sum k_i$. 

Since $x^{k_i}$ is  not an SDM, support of $\HF_{*}(H^{\nat k_i}, \bx^{k_i})$ is contained in the half open interval $[k_i\Delta(\bx)-n, k_i\Delta(\bx)+n) $. So there exists $\delta_{k_i} > 0$ depending on $k_i$ but not on the capping, such that $l_i -k_i\Delta(\bx) -n \leq -\delta_i$. Summing over $i$ gives 
$$\sum_{i=1}^r l_i -k\Delta(\bx) -rn = l-k\Delta(\bx) -n \leq-\sum_{i=1}^r \delta_{k_i} \leq -r\delta,
$$ 
where $\delta =\min\{ \delta _{k_i} \, \vert \, k_i \leq N\} $. So when $ r> 2n / \delta $, $l$ goes out of the support of $\HF_*(H^{\nat k}, \bx^k)$.

\begin{Remark}
When the mean index is an integer we can choose $\delta=1$ independent of $N$. Thus, in that case the algebra is uniformly nilpotent. But this is not true in general. For an example, consider the autonomous Hamiltonian $H(x,y)=-\lambda (x^2+y^2)$ where $\lambda >0$ is a small irrational number.  For $C^2$-small autonomous Hamiltonians, e.g. $H$, pair-of-pants product in the local Floer homology of an isolated critical point agrees with the dual of cup-product in the local Morse  homology of the same point. The latter is non-zero for a local maximum. So the map
\[
\HF_*(H,0) \otimes \cdots \otimes \HF_*(H,0) \rightarrow  \HF_*(H^{\nat k},0)
\]
given by $k$-fold pair-of-pants product is non-zero provided that $H^{\nat k}$ is $C^2$-small. Here notion of being $C^2$-small changes with iteration, since we are considering $H^{\nat k}$ as $k$-periodic (since we are not normalizing the period).  Now let $l \in \N$ such that $ \lambda \gg  l\lambda-\lfloor l\lambda \rfloor =\theta > 0$, and let $K(x,y)=-\theta (x^2 +y^2)$. The isomorphism 
\[
 \HF_*(K,0)  \rightarrow \HF_*(H^{\nat l}, 0)
\]
given by composing with the loop diffeomorphism generated by the Hamiltonian $G(x,y)= -\lfloor l\lambda \rfloor (x^2+y^2)$ preserves the pair-of-pants product (see Section 5.2 in \cite{GG:gap}). Since the iterated Hamiltonian $K^{\nat k}$ stays $C^2$-small for higher iterations compare to $H^{\nat k}$, the map
 \[
\HF_*(H^{\nat l},0) \otimes \cdots \otimes \HF_*(H^{\nat l},0) \rightarrow  \HF_*(H^{\nat lk},0)
\]
would not vanish for larger $k$. Hence the algebra $\FA(H,0)$ is not uniformly nilpotent. 
\end{Remark}

\begin{Remark}
Conclusion of Theorem \ref{thm:2} is false in the presence of a symplectically degenerate maximum. In fact, see Section 5.2 in \cite{GG:gap} for details, if $w \in \HF_{\Delta(\bx)+n}(H,\bx)$ is the generator of the top degree homology (which is one-dimensional) of an SDM orbit $\bx$, then for every large prime $k$ the $k$th power of $w$ is non-zero in  $\HF_{\Delta(\bx^k)+n}(H^{\nat k},\bx^k)$.
\end{Remark}

\begin{Remark}
Roughly speaking, the only restriction on the local Morse homology of an isolated critical point is that with all the algebraic structures it is isomorphic to the homology of a suspension, see \cite{CLOT, Pe}. In particular, the cup-product and the Massey products vanish in the local Morse (co)homology, \cite{Vi}, but other cohomology operations need not be trivial (e.g. Steenrod squares). Non-uniform nilpotency of the pair-of-pants product is a Floer theoretic counterpart of this vanishing phenomena. However, similar to the local Morse homology, the local Floer homology can carry many non-nilpotent cohomology operations.
\end{Remark}

\subsection{Proof of Theorem \ref{thm:1}}
\label{sec:pf1} 

As discussed in Section \ref{sec:notation}, the mean action spectrum $\widehat{\CS}(H)$ of $H \colon S^1 \times M \rightarrow \reals$ is infinite when $\omega$ is not aspherical; or in the existence of an SDM orbit, see  \cite{Gi:CC} for details. In this section, we assume that $\omega$ is aspherical and none of the orbits of $\varphi_H$ is an SDM. 

Arguing by contradiction, assuming in addition that $\widehat{\CS}(H)$ is finite, let $\widehat{\CS}(H)= \{a_1, \dots , a_m \}$ be the mean action spectrum (ordered) of $H$. Choose $c_i \in \reals$ such that $a_i < c_i < a_{i+1}$. By the energy estimates for the product, see Section \ref{sec:product} and \cite{GG:conley} for details, the filtered Floer algebra  $\FA^{(-\infty,\, c_1)}(H)$ splits
\[
\FA^{(-\infty,\, c_1)}(H)=\bigoplus_{\Aa_{H^{\nat k}}(x_i)=ka_1} \FA(H^{\nat k}, x_i)
\]
as an algebra. Each summand is non-uniformly nilpotent by Theorem \ref{thm:2} and there are finitely many summands less than a fixed iteration, so the sum is non-uniformly nilpotent. Next we will inductively argue that the total algebra $\FA(H)$ is non-uniformly nilpotent (or just nilpotent, since summands are isomorphic), which is impossible (see Section \ref{sec:product}).

Recall that a graded algebra $\bigoplus_{k=1}^{\infty} \mathcal{A}_k$ is called non-uniformly nilpotent, if for all $N\in \N$ there is $m\in \N$ such that all $m$-fold products vanish when restricted to $\bigoplus_{k=1}^{N} \mathcal{A}_k$. We call the minimum such $m$ the \emph{degree of nilpotentcy} of $\bigoplus_{k=1}^{N} \mathcal{A}_k$. Fix $N \in \mathbb{N}$ and let $l$, $k$ be the degrees of nilpotency of 
\[
\FA^{(c_i,\, c_{i+1})}(H):=\bigoplus_{\Aa_{H^{\nat k}}(x_j)=ka_{i+1}} \FA(H^{\nat k} x_j)
\]
and $\FA^{(-\infty, c_{i})}(H)$ when restricted to first $N$, $Nl$ iterations respectively. We will show that $\FA^{(-\infty, c_{i+1})}(H)$ is nilpotent with degree at most $k l$ when restricted to first $N$ iterations, and hence non-uniformly nilpotent. 

Let $A,B,C$ be the Floer chain complexes which give rise to Floer algebras above, i.e., $\H_{*}(A)=\FA^{(-\infty, c_{i})}(H)$, $\H_{*}(B)=\FA^{(c_i, c_{i+1})}(H)$, $\H_{*}(C)=\FA^{(-\infty, c_{i+1})}(H)$. Let $D$ be the quotient chain complex coming from the inclusion $A \otimes A \rightarrow C \otimes C$. Using the auxilary complex $D$ we form a commutative diagram
\[
\begin{tikzcd}
0 \arrow{r}{} & A \otimes A \arrow{r}{} \arrow{d}{} & C \otimes C \arrow{r}{} \arrow{d}{} & \arrow {d}{} D \arrow{r}{} & 0 \\
0 \arrow{r}{} & A \arrow{r}{} & C \arrow{r}{} & B \arrow{r}{} & 0
\end{tikzcd}
\]
with exact rows, where the vertical arrows are pair-of-pants product. It induces a commutative diagram
\[
\begin{tikzcd}
\HF_{*}(A) \otimes  \HF_{*}(A) \arrow{r}{} \arrow{d}{} &  \HF_{*}(C) \otimes  \HF_{*}(C) \arrow{r}{} \arrow{d}{} & \arrow {d}{}  \HF_{*}(D) \\
\HF_{*}(A) \arrow{r}{P} &  \HF_{*}(C) \arrow{r}{R} &  \HF_{*}(B) 
\end{tikzcd}
\]
in homology. Note that the product map $D \rightarrow B$ factors through $B\otimes B$; so $\HF_{*}(D)$ has the same product structure with $\HF_{*}(B) \otimes  \HF_{*}(B)$.

Now take $k l$ classes $w_i \in \HF_{*}(C)$ that belong to first $N$ iterations. If at least $k$ of them are in the image of $P$, using the associativity of the product and the commutativity of the first block in the diagram, we conclude that $w_1\cdot \dots {\cdot w_{kl}}=0$. 

If not, take $l$ of the classes $w_{i_j}$ that are not in the image of $P$. This time using commutativity of the second block, we conclude that the product $w_{i_1}\cdot \dots {\cdot w_{i_l}}$ is in the kernel of $R$ and hence in the image of $P$. There are at least $k-1$ such $l$-tuples; and if that is case, then there is at least one class among the remaining $l$-classes, that is in the image of $P$. So in total we obtain $k$ classes (from first $Nl$ iterations) in the image of $P$ and we go back to the first case to conclude that  $w_1\cdot \dots {\cdot w_{kl}}=0$. 

\begin{Remark}
In fact a slightly stronger result holds. Let $c \in \R$ such that  $c \notin \widehat{\CS}(H)$ and $c$ bounds from above the action spectrum of $H$. Consider the commutative diagram
\[
\begin{tikzcd}
\HF_*(H) \otimes \cdots \otimes \HF_*(H)   \arrow{r}{} \arrow{d}{} & \HF^{(-\infty, \, kc)}_*(H^{\nat k}) \arrow{d}{}  \\
\HF_*(H) \otimes \cdots \otimes \HF_*(H)  \arrow{r}{} &  \HF_*(H^{\nat k})
\end{tikzcd}
\]
where horizontal arrows are $k$-fold pair-of-pants product, the left vertical arrow is the identity map, and the right vertical arrow is the induced map coming from the inclusion of the subcomplex  $\CF^{(-\infty, \, kc)}_*(H^{\nat k})$. If $\widehat{\CS}(H)$ has no accumulation point, then using the same argument as above, one can show that the first row in the diagram is identically zero for sufficiently large $k$. Then the second row is also zero, which is not possible since quantum cup-product has a unit. 
\end{Remark}

\end{document}